\documentclass[10pt]{amsart} 
\usepackage{amscd}  
\usepackage{amsfonts}  
\usepackage{amssymb}  
\usepackage{latexsym}

\newcommand{\ncm}{\newcommand} 
 
\ncm{\beq}{\begin{equation}} 
\ncm{\eeq}{\end{equation}} 
 
\newtheorem{thm}{Theorem}[section] 
 
\newtheorem{lem}[thm]{Lemma} 
\newtheorem{cor}[thm]{Corollary} 
\newtheorem{lem&def}[thm]{Lemma \& Definition} 

\theoremstyle{definition}
\newtheorem{defi}[thm]{Definition}

\theoremstyle{remark}

\def\Set{\mathsf{Set}} 
 
\def\Ab{\mathsf{Ab}}

\def\Mon{\mathsf{Mon}}

\def\Cat{\mathsf{Cat}} 
\def\MonCat{\mathsf{MonCat}}
\def\OpmonCat{\mathsf{Mon_{op}Cat}}
\def\Bmd{\mathsf{Bmd}}
\def\Bgd{\mathsf{Bgd}}

\def\Fun{\mathsf{Fun}}

\def\M{\mathcal{M}} 
 
\def\V{\mathcal{V}}

\def\T{\mathcal{T}} 
\def\C{\mathcal{C}} 

\def\ttT{\mathtt{T}}
\def\ttL{\mathtt{L}}
\def\ttR{\mathtt{R}}
\def\ttM{\mathtt{M}}

\ncm{\Ann}{\mbox{\rm Ann}} 
\ncm{\End}{\mbox{\rm End}\,} 
 
\def\Hom{\mbox{\rm Hom}\,}

\def\to{\rightarrow} 
 
\def\b{\,\Box\,}    
\def\c{\circ}
\def\o{\otimes}
\def\x{\times}                 

\def\bra{\langle} 
\def\ket{\rangle} 
\def\under{\mbox{\rm\_}\,} 
\ncm{\rarr}[1]{\stackrel{#1}{\longrightarrow}} 
\ncm{\larr}[1]{\stackrel{#1}{\longleftarrow}} 
\def\iso{\rarr{\sim}}
\def\cop{\Delta}

\def\eps{\varepsilon} 
 
\def\du1{\hat 1}

\def\0{_{(0)}} 
\def\1{_{(1)}} 
\def\2{_{(2)}} 
\def\3{_{(3)}}

\def\du1{\hat 1} 
 
\def\ract{\triangleleft} 
\def\op{^{\rm op}}
\def\a{\mathbf{a}}
\def\l{\mathbf{l}}
\def\r{\mathbf{r}}

\renewcommand{\theequation}{\thesection.\arabic{equation}} 

\begin{document} 
 
\title{Adjointable monoidal functors and quantum groupoids}
 
\author{Korn\'el Szlach\'anyi}  
\address{Theory Division, Research Institute for Particle and Nuclear Physics,
Budapest, H-1525 Budapest, P. O. Box 49, Hungary}  

\begin{abstract}
Every monoidal functor $G\colon\C\to\M$ has a canonical factorization through
the category $_R\M_R$ of bimodules in $\M$ over some monoid $R$ in $\M$ in
which the factor $U\colon\C\to\,_R\M_R$ is strongly unital. 
Using this result and the characterization of the forgetful functors
$\M_A\to\,_R\M_R$ of bialgebroids $A$ over $R$ given by Schauenburg
\cite{Schauenburg: Bial} together with their bimonad description given by the
author in \cite{Sz: JPAA} here we characterize the "long" forgetful functors
$\M_A\to\,_R\M_R\to\M$ of both bialgebroids and weak bialgebras.  
\end{abstract}

\thanks{
Based on a talk given at the conference
"Hopf Algebras in Noncommutative Geometry and Physics" organized by Stefaan
Caenapeel and Fred van Oystaeyen, Brussels, May 28 - June 1, 2002.
}

\maketitle

\section{Introduction}

Takeuchi's $\x_R$-bialgebras \cite{Takeuchi} or, what is the same
\cite{Brz-Mil}, Lu's bialgebroids \cite{Lu} provide far reaching
generalizations of the notion of bialgebra. A bialgebroid $A$ is, roughly
speaking, a bialgebra over some non-commutative $k$-algebra $R$. 
With noncommutativity of $R$, however, a new phenomenon appears: the separation
of algebra and coalgebra structures into two different categories.
While bialgebras are monoids and comonoids in the same category
$\M_k$, bialgebroids are monoids in $\M_k$ (or in $_{R^e}\M_{R^e}$) but
comonoids in $_R\M_R$. This makes the compatibility conditions rather
difficult to formulate. 

Simplification, if at all, is expected on passing to the level of categories
and functors. Let us take, for example, a $k$-algebra $A$ and associate to it
the monad $T=\under\o A$ on $\M_k$. The monads obtained that way are precisely
the monads which have right adjoints. For this characterization of algebras the
closed monoidal structure of $\M_k$ is essential. The Eilenberg-Moore category
of "$T$-algebras" \cite{MacLane,TTT}, or perhaps better to say, "$T$-modules"
is nothing but the category $\M_A$ of right $A$-modules equipped with the
forgetful functor $\M_A\to\M_k$. Similarly, we can consider monads on the
closed monoidal category $_R\M_R\equiv\M_{R^e}$, where $R^e=R\op\o R$, with
right adjoints. It turns out that these are, up to isomorphisms, precisely the
monads $\under\o_{R^e}A$ associated to monoids $A$ in $_{R^e}\M_{R^e}$, also
called $R^e$-rings. The forgetful functor $U^A\colon\M_A\to\,_R\M_R$ is
monadic, has a left adjoint but has no monoidal structure. At this
point bialgebroids enter naturally via Schauenburg's theorem
\cite{Schauenburg: Bial}: the monoidal structures on $\M_A$ such that $U^A$ is
strict monoidal are in one-to-one correspondence with (right) bialgebroid
structures on the $R^e$-ring $A$.

Monoidal structures on $\M_A$ can also be described by
opmonoidal structures on the monad $T=\under\o_{R^e}A$. In a recent paper
\cite{Sz: JPAA} bialgebroids have been characterized as the \textit{bimonads}
on $_R\M_R$ the underlying functors of which have right adjoints. A bimonad, or
opmonoidal monad \cite{Moerdijk,McCrudden}, is a monad in the 2-category
$\Mon_{op}\Cat$ of monoidal categories, opmonoidal functors and opmonoidal
natural transformations. More explicitly, a bimonad $\bra
T,\gamma,\pi,\mu,\eta\ket$ on a monoidal category $\bra \M,\o,i\ket$ consists
of  
\begin{enumerate}
\item an endofunctor $T\colon\M\to\M$
\item a natural transformation $\gamma_{x,y}\colon T(x\o y)\to Tx\o Ty$
\item an arrow $\pi\colon Ti\to i$
\item a natural transformation $\mu_x\colon T^2x\to Tx$
\item and a natural transformation $\eta_x\colon x\to Tx$
\end{enumerate}
such that $\bra T,\mu,\eta\ket$ is a monad, $\bra T,\gamma,\pi\ket$ is an
opmonoidal functor, i.e., a monoidal functor in $\bra \M\op,\o,i\ket$, and
$\mu$ and $\eta$ are opmonoidal natural transformations in the obvious sense.
Bimonads, and therefore bialgebroids, too, form a 2-category $\Bmd$
and $\Bgd\subset\Bmd$, respectively \cite{Sz: JPAA}. 

The forgetful functors $U^A\colon\M_A\to\,_R\M_R$ of bialgebroids over $R$ can
be characterized as the strong monoidal monadic functors to $_R\M_R$ that have
right adjoints \cite[Corollary 4.16]{Sz: JPAA}.
In this paper we will study analogue characterizations of the \textit{long
forgetful functors} $G^A\colon\M_A\to\M_k$. This is motivated by situations
where the base algebra $R$ is not given a priori. It is also closer to the
classical Tannaka-Krein situation where one reconstructs the "grouplike"
object $A$ as the set of natural transformations $G^A\to G^A$. Apart from set
theoretical controversies ($\M_A$ is not small, which will be compensated by
assuming the existence of left adjoints for our functors) this reconstruction
is possible for the long forgetful functor $G^A$ but not for the short
forgetful functor $U^A$. Of course, the novelty is that now the long forgetful
functor is not strong monoidal.

To recover $R$ from $G^A$ is in fact very easy. One takes the image of
the unit object (the trivial $A$-module) under $G^A$. Since the unit object is
always a monoid, it is mapped by the monoidal forgetful functor to a monoid in
$\M_k$. This gives us the algebra $R$. This construction is possible for any
monoidal functor and a closer look will show in Section \ref{sec: can fac} that,
under mild assumptions on $\M$, every monoidal functor $G:\C\to\M$ can be
factorized as $\C\rarr{U}\,_R\M_R\to\M$ with $U$ monoidal but strictly unital. 

The monoidal functors $G$ for which $U$ is strong monoidal will be called 
\textit{essentially strong monoidal}. Clearly, the $G^A$ of a bialgebroid is
an example of such functors. Finding the extra conditions on an essentially
strong monoidal functor $G\colon\C\to\M$ that makes it (1) either factorize
through the long forgetful functor $G^A$ of a unique bialgebroid (2) or become
isomorphic to such a $G^A$ is part of a Tannaka duality program for
bialgebroids. This has been carried out for "short" forgetful functors in 
\cite{Sz: JPAA}. This type of duality theory uses monad theory to 
characterize the large module categories of quantum groupoids together with
their forgetful functors. With the results of the present paper we make some
small steps in the direction of extending Tannaka theory from strong monoidal
to monoidal functors. As for the state of the art of the traditional method we
have to mention the recent papers by Ph\`ung H{\`o} H\'ai \cite{Phung Ho Hai}
and another one by Hayashi \cite{Hayashi: Tannaka} which prove Tannaka duality
theorems for Hopf algebroids and for face algebras, respectively. In their
approach, as in that of Saavedra-Rivano, Deligne, Ulbrich and others (see
\cite{Pareigis: NCAG,Joyal-Street: Tannaka}) small categories are equipped
with strong monoidal functors to a (sometimes rigid) category of bimodules and
the task is to find a universal factorization through the \textit{comodule}
category of a quantum groupoid. 

\newpage

The organization of the paper is as follows. In Section \ref{sec: can fac} 
we pove the canonical factorization of general monoidal functors through a
bimodule category. After
touching the general case of long forgetful functors of bimonads in Section
\ref{sec: EM} we determine a class of essentially strong monoidal functors in
Section \ref{sec: bgd} which factorize through the $G^A$ of a bialgebroid.
Although we present the proof over the base category $\M_k$, $\Ab$ or $\Set$,
it is often indicated how it could be extended to a general base category
$\V$. In this way we intend to make playground for exotic examples of
bialgebroids. Then in Section \ref{sec: G^A} the long forgetful functors of
bialgebroids are characterized up to equivalence. Finally, in Section
\ref{sec: WBA}, the special case of weak bialgebras are considered, now over
$\M_k$, the long forgetful functors of which can be recognized as those that
have both monoidal and opmonoidal structures and these two obey compatibility
conditions that can be called a \textit{separable Frobenius structure} on the
forgetful functor. This characterization of weak bialgebra forgetful functors
was already sketched in \cite{Sz: Siena} calling them "split monoidal"
functors.

\section{The canonical factorization of monoidal functors} \label{sec: can fac}

Let $\C$ be a monoidal category with monoidal product $\b\colon\C\x\C\to\C$
and unit object $e\in\C$. Then we have coherent natural isomorphisms
$\a_{a,b,c}\colon a\b(b\b c)\iso(a\b b)\b c$, $\l_c\colon e\b c\iso c$ and
$\r_c\colon c\b e\iso c$ satisfying $\l_e=\r_e$, the triangle and the pentagon
identity. A monoid $\bra m,\mu,\eta\ket$ in $\C$ is an object $m$ together
with arrows $\mu\colon m\b m\to m$, $\eta\colon e\to m$ satisfying
associativity and unit axioms. There is always a canonical monoid: the unit
object $e$ equipped with multiplication 
$l_e\colon e\b e\to e$ and unit the identity arrow $e\colon e\to e$.
Moreover, every object $c$ of $\C$ is a bimodule over this canonical monoid
via the actions $\l_c\colon e\b c\to c$ and $\r_c\colon c\b e\to c$. The
bimodule axioms follow simply from recognizing that the three associativity
axioms  \begin{align}
\l_c\c(e\b\l_c)&=\l_c\c(\l_e\b c)\c\a_{e,e,c}\\
\l_c\c(e\b\r_c)&=\r_c\c(\l_c\b e)\c\a_{e,c,e}\\
\r_c\c(\r_c\b e)&=\r_c\c(c\b\l_e)\c\a_{c,e,e}^{-1}
\end{align}
are consequences of special cases of the 
triangle diagrams valid in any monoidal category while the unit axioms become
identities.
In this sense every monoidal category is a category of bimodules. The more
precise statement will be clear after applying the Theorem below to the
identity functor of $\C$. Let us recall an important property of
monoidal functors: They map monoids to monoids and (bi)modules to (bi)modules.

\begin{lem} \label{lem: mon to mon}
Let $\bra G,G_2,G_0\ket\colon\bra\C,\b,e\ket\to\bra\M,\o,i\ket$ be
a monoidal functor. 
\begin{enumerate}
\item If $\bra m,\mu,\eta\ket$ in $\C$ is a monoid in $\C$ then
$\bra Gm,G\mu\c G_{m,m},G\eta\c G_0\ket$ is a monoid in $\M$. 
\item Let $m$ and $n$ be monoids in $\C$. If $\bra b,\lambda,\rho\ket$
is an $m$-$n$ bimodule in $\C$ then the triple $\bra Gb,G\lambda\c
G_{m,b},G\rho\c G_{b,n}\ket$ is a $Gm$-$Gn$ bimodule in $\M$.
\end{enumerate}
\end{lem}
\begin{proof}
(1) is well-known and can be found e.g. in \cite{Street: Quantum Groups}.
(2) is also known to many authors although an explicit proof is difficult
to find. Just to advertise the statement we compute here commutativity of
the left and right actions: 
\begin{align*}
\lambda'\c(Gm\o\rho')&=G\lambda\c G_{m,b}\c(Gm\o G\rho)\c(Gm\o G_{b,n})\\
&=G\lambda\c G(m\b\rho)\c G_{m,b\b n}\c(Gm\o G_{b,n})\\
&=G\rho\c G(\lambda\b n)\c G\a_{m,b,n}\c G_{m,b\b n}\c(Gm\o G_{b,n})\\
&=G\rho\c G(\lambda\b n)\c G_{m\b b,n}\c(G_{m,b}\o Gn)\c\a_{Gm,Gb,Gn}\\
&=G\rho\c G_{b,n}\c (G\lambda\o Gn)\c(G_{m,b}\o Gn)\c\a_{Gm,Gb,Gn}\\
&=\rho'\c(\lambda'\o Gn)\c\a_{Gm,Gb,Gn}
\end{align*}
\end{proof}
Dually, opmonoidal functors map comonoids to comonoids and (bi)comodules
to (bi)comodules.

Sofar $\M$ was an arbitrary monoidal category. In order for the category 
$_m\M_m$ of bimodules in $\M$ over a monoid $m$ in $\M$ to have a monoidal
structure we need the assumptions that $\M$ has coequalizers and the tensor
product $\o$ preserves coequalizers in both arguments. This is because the
tensor product over $m$ of a right $m$-module $\rho_x\colon x\o m\to x$
with a left $m$-module $\lambda_y\colon m\o y\to y$ is a coequalizer
\beq\label{dia: tensorproduct} 
\parbox[c]{4.5in}
{
\begin{picture}(320,65)(-20,-10) 
\put(0,40){$x\o (m \o y)$}
\put(0,0){$(x\o m)\o y$}
\put(30,35){\vector(0,-1){25}}
\put(33,22){$\wr$}
\put(115,20){$x\o y$}
\put(65,40){\vector(3,-1){45}}
\put(65,5){\vector(3,1){45}}
\put(78,38){$\scriptstyle x\o\lambda_y$}
\put(78,4){$\scriptstyle \rho_x\o y$}
\put(150,24){\line(1,0){45}} \put(194,21.4){$\twoheadrightarrow$}
\put(210,20){$x\o_m y$}
\end{picture}
}
\eeq
The construction of the monoidal product $\o_m$ on $_m\M_m$ together with
coherence isomorphisms is a long but standard procedure. At the
end one obtains a monoidal category $\bra_m\M_m,\o_m,m\ket$ together with a
monoidal forgetful functor $\Gamma^m\colon\,_m\M_m\to\M$ sending the bimodule
$\bra x,\lambda,\rho\ket$ to its underlying object $x$. The monoidal
structure 
\begin{align}
\Gamma^m_{x,y}\colon&\Gamma^m(x)\o \Gamma^m(y)\to\Gamma^m(x\o_m y)\\
\Gamma^m_0\colon& i\to \Gamma^m(m)
\end{align}
is provided by the chosen coequalizers and by the unit $i\to m$ of the
monoid $m$, respectively. 

If $\sigma\colon m\to n$ is a monoid morphism then there is a functor
$\Gamma^\sigma\colon\,_n\M_n\to\,_m\M_m$ mapping the bimodule $\bra
x,\lambda,\rho\ket$ to the bimodule $\bra x,\lambda\c(\sigma\o
x),\rho\c(x\o\sigma)\ket$. So, $\Gamma^\sigma\Gamma^m=\Gamma^n$.
This defines a functor $\Gamma\colon\Mon\M\to\MonCat/\M$.

\begin{thm}\label{thm: monfunc}
Let $\bra \M,\o,i\ket$ be a monoidal category with coqualizers
and such that $\o$ preserves coequalizers in both arguments. If
$\bra \C,\b,e\ket$ is a monoidal category and
$G\colon\C\to\M$ is a monoidal functor
then there is a monoid $\bra R,\mu^R,\eta^R\ket$ in $\M$ and a strictly unital
monoidal functor $U\colon\C\to\,_R\M_R$ such that 
\begin{enumerate}
\item $G$, as a monoidal functor, can be factorized as $\Gamma^RU$, i.e.,
\begin{align}
G&=\Gamma^RU \label{eq: monfac1}\\
G_{a,b}&=\Gamma^R U_{a,b}\c \Gamma^R_{Ua,Ub} \label{eq: monfac2}\\
G_0&=\Gamma^R U_0\c \Gamma^R_0 \label{eq: monfac3}
\end{align}
\item If $S$ is a monoid in $\M$ and $V\colon\C\to\,_S\M_S$ is a monoidal
functor such that $G=\Gamma^SV$, as monoidal functors, then there exists a unique
monoid morphism $\sigma\colon S\to R$ such that $V=\Gamma^\sigma U$.
\end{enumerate}
\end{thm}
\begin{proof}
By Lemma \ref{lem: mon to mon} the image under $G$ of the unit monoid $e$ is
a monoid $R=\bra Ge,\mu^R,\eta^R\ket$ with underlying object $Ge$. Also by the
Lemma, every object $c$ in $\C$, as an $e$-$e$-bimodule, is mapped by $G$ to
the $R$-$R$-bimodule  
\beq
Uc:=\bra Gc,\lambda_{Uc},\rho_{Uc}\ket
\eeq
where
\begin{align}
\label{eq: lambda}
\lambda_{Uc}:=\quad\begin{CD}Ge\o Gc@>G_{e,c}>>G(e\b c)@>G\l_c>>Gc\end{CD}\\
\label{eq: rho}
\rho_{Uc}:=\quad\begin{CD}Gc\o Ge@>G_{c,e}>>G(c\b e)@>G\r_c>>Gc\end{CD}
\end{align}
Since these actions are natural in $c$, the $G\tau$ of every arrow $\tau\colon
c\to d$ lifts to a bimodule morphism $U\tau\colon Uc\to Ud$. This defines the
functor $U$ which obviously satisfies (\ref{eq: monfac1}). In order to define a
monoidal structure for $U$ notice that 
\beq
G_{a,b}\c(\rho_{Ua}\o Gb)\c\a_{Ga,Ge,Gb}=G_{a,b}\c(Ga\o\lambda_{Ub})
\eeq
holds true as a consequence of the hexagon of $G_2$. Therefore, universality
of the coequalizer implies the existence of a unique arrow $Ga\o_R Gb\to
G(a\b b)$ such that
\beq
\begin{picture}(200,100)(0,-5)
\put(0,70){$Ga\o Gb$}
\put(40,72){\vector(1,0){80}} \put(70,78){$\Gamma^R_{Ua,Ub}$}
\put(125,70){$Ga\o_R Gb\equiv\Gamma^R(Ua\o_R Ub)$}
\multiput(145,62)(0,-4){13}{\line(0,1){2}}
\put(145,9){\vector(0,-1){0}} \put(150,32){$\Gamma^R U_{a,b}$}
\put(20,65){\vector(3,-2){98}} \put(48,27){$G_{a,b}$}
\put(125,-5){$G(a\b b)$}
\end{picture}
\eeq
commutes. Moreover, this new arrow lifts to a bimodule morphism $U_{a,b}$ since
the other two in this diagram also lift to $_R\M_R$ and $\o$ preserves
coequalizers. Setting $U_0\colon \,_RR_R\to Ue$ to be the identity arrow we
obtain a monoidal functor $U\colon\C\to\,_R\M_R$ for which (\ref{eq: monfac2})
and (\ref{eq: monfac3}) hold and the unit of which, $U_0$, is an identity
arrow, i.e., it is strictly unital. This proves property (1).  
To prove the universal property (2) we start with uniqueness. If $\sigma$
exists such that $\Gamma^\sigma U=V$ is a factorization in $\MonCat$ then, in
particular, $V_0=\Gamma^\sigma U_0\c\Gamma^\sigma_0$. Since $U$ is strict
unital, $\Gamma^SV_0=\Gamma^S\Gamma^\sigma_0$. This means that the underlying
arrow in $\M$ of the bimodule morphism $\Gamma^\sigma_0\colon
S\to\,_{\sigma(S)}R_{\sigma(S)}$ is uniquely determined by that of $V_0$. But
this bimodule morphism is the same as the monoid morphism
$\sigma\colon S\to R$ (as arrows of $\M$). To prove existence we therefore
define $\sigma:=\Gamma^SV_0\colon
\Gamma^SS\to\Gamma^RR$, for the time being as an arrow in $\M$. It preserves
the unit due to the monoidal factorization, 
\[
\sigma\c\eta^S=\Gamma^SV_0\c\Gamma^S_0=G_0=\Gamma^R_0=\eta^R
\]
and it is multiplicative because 
\begin{align*}
\mu^R\c(\sigma\o\sigma)&=G\l_e\c G_{e,e}\c(\Gamma^SV_0\o\Gamma^SV_0)\\
&=\Gamma^SV\l_e\c\Gamma^SV_{e,e}\c\Gamma^S_{Ve,Ve}\c(\Gamma^SV_0\o\Gamma^SV_0)\\
&=\Gamma^SV\l_e\c\Gamma^SV_{e,e}\c\Gamma^S(V_0\o_S V_0)\c\Gamma^S_{S,S}\\
&=\Gamma^S\l_{Ve}\c\Gamma^S(S\o_S V_0)\c\Gamma^S_{S,S}
=\Gamma^SV_0\c\Gamma^S\l_S\c\Gamma^S_{S,S}\\
&=\sigma\c\mu^S\,.
\end{align*}
So $\sigma$ lifts to a monoid morphism $S\to R$. The proof of $\Gamma^\sigma
U=V$ requires to show that $\lambda_{Uc}\c(\sigma\o Gc)=\lambda_{Vc}$ and
$\rho_{Uc}\c(Gc\o \sigma)=\rho_{Vc}$ for all object $c$ in $\C$. E.g.,
\begin{align*}
\lambda_{Vc}&=\Gamma^S\l_{Vc}\c\Gamma^S_{S,Vc}\\
&=\Gamma^SV\l_c\c\Gamma^SV_{e,c}\c\Gamma^S(V_0\o_S V_0)\c\Gamma^S_{S,Vc}\\
&=G\l_c\c\Gamma^SV_{e,c}\c\Gamma^S_{Ve,Vc}\c(\Gamma^SV_0\o Gc)\\
&=G\l_c\c G_{e,c}\c(\Gamma^SV_0\o Gc)\\
&=\lambda_{Uc}\c(\sigma\o Gc)\,.
\end{align*}
In order to have $\Gamma^\sigma U=V$ in $\MonCat$ we still have to show
$V_{a,b}=\Gamma^\sigma U_{a,b}\c\Gamma^\sigma_{Ua,Ub}$ and
$V_0=\Gamma^\sigma U_0\c\Gamma^\sigma_0$. The latter follows directly from the
definition of $\sigma$ while the former can be shown by applying the
faithful $\Gamma^S$ on both hand sides and then composing with the coequalizer
$\Gamma^S_{Va,Vb}$ from the right. Since a coequalizer is epi, it is right
cancellable and the statement follows from (\ref{eq: monfac2}).
\end{proof}
We call $G=\Gamma U$ the canonical factorization of the monoidal functor $G$.

\begin{defi} \label{def: ess}
A monoidal functor $G$ is called essentially strong monoidal if the $U$ in its
canonical factorization $G=\Gamma U$ is strong monoidal.   
\end{defi}

Equivalently, $G$ is essentially strong monoidal if for all objects $a,b$ of
$\C$ the arrow $G_{a,b}$ is a coequalizer in the diagram
\beq\label{exact}
\parbox[c]{4.5in}{
\begin{picture}(320,65)(-20,-10) 
\put(0,40){$Ga\o (Ge \o Gb)$}
\put(0,0){$(Ga\o Ge)\o Gb$}
\put(40,35){\vector(0,-1){25}}
\put(43,22){$\wr$}
\put(135,22){$Ga\o Gb$}
\put(85,40){\vector(3,-1){45}}
\put(85,5){\vector(3,1){45}}
\put(98,38){$\scriptstyle Ga\o \lambda_{Ub}$}
\put(98,5){$\scriptstyle \rho_{Ua}\o Gb$}
\put(192,25){\vector(1,0){40}}
\put(200,31){$G_{a,b}$}
\put(240,22){$G(a\b b)$}
\end{picture}
}
\eeq
where $\lambda_{Uc}$, $\rho_{Uc}$ denote the left and right $R$-actions
(\ref{eq: lambda}) and (\ref{eq: rho}) on $Gc$.

\section{The relation with the Eilenberg-Moore construction} \label{sec: EM}

Let us recall some basic facts about monad theory \cite{MacLane,TTT}.
Any functor $G\colon \C\to\M$ with a left adjoint $F$ determines a monad
$\ttM =\bra GF,G\eps F,\eta\ket$ on $\M$ where $\eps\colon FG\to\C$ is the counit
and $\eta\colon\M\to GF$ is the unit of the adjunction. The Eilenberg-Moore
construction associates to any monad $\ttM $ on $\M$ a category $\M^\ttM $ 
the objects of which are the
$\ttM $-algebras $\bra x,\alpha\ket$ where $x$ is an object of
$\M$ and $\alpha\colon \ttM x\to x$ is an arrow of $\M$ satisfying $\alpha\c
M \alpha=\alpha\c\mu_x$ and $\alpha\c\eta_x=x$. The arrows of $\M^\ttM $
from $\bra x,\alpha\ket$ to $\bra y,\beta,\ket$ are the arrows $\tau\colon
x\to y$ in $\M$ for which $\tau\c\alpha=\beta\c M \tau$.

The forgetful functor $U^\ttM \colon\M^\ttM \to\M$ sending $\bra x,\alpha\ket$ to $x$
has a left adjoint such that the monad associated to this adjunction is
precisely the original $\ttM $. 

Any functor $G\colon \C\to\M$ with a left adjoint can be factorized as $G=U^\ttM 
K$, where $\ttM $ is the monad of the adjunction, with a unique
$K\colon\C\to\M^\ttM $, called the comparison functor. Explicitly,
\beq
c\stackrel{K}{\longmapsto}\bra Gc,G\eps_c\ket\stackrel{U^\ttM }{\longmapsto}Gc\,.
\eeq
If $K$ is an equivalence of categories the $G$ is called monadic.

Now assume that we have a monoidal functor $G\colon\C\to\M$ with the underlying
ordinary functor having a left adjoint $F$. We briefly say that $G$ is a right
adjoint monoidal functor. Then we have two factorizations of $G$, the
$G=\Gamma U$ provided by Theorem \ref{thm: monfunc} and $G=U^\ttM  K$ provided by
the Eilenberg-Moore construction. In this situation one expects a relation
between $\ttM $-algebras and $R$-$R$-bimodules.

Every monoid in $\M$, therefore $R$ too, provides two monads $\ttL=R\o\under$
and $\ttR=\under\o R$. We can define two monad morphisms
\begin{align}
\lambda^\ttM \colon& \ttL\to \ttM \qquad
\lambda^\ttM _x=G(\l_{Fx})\c G_{e,Fx}\c(R\o\eta_x)\\
\rho^\ttM \colon& \ttR\to \ttM \qquad
\rho^\ttM _x=G(\r_{Fx})\c G_{Fx,e}\c(\eta_x\o R) 
\end{align}

They commute in the sense of the diagram
\beq
\begin{CD}
\ttL\ttR@>\lambda^\ttM \rho^\ttM >>\ttM ^2@>G\eps F>>\ttM \\
@VV\wr V @. @|\\
\ttR\ttL@>\rho^\ttM \lambda^\ttM >>\ttM ^2@>G\eps F>>\ttM 
\end{CD}
\eeq
where the vertical isomorphism can be obtained from the associator as
$\a_{R,\under,R}$.

\begin{lem}
If $\bra x,\alpha\ket$ is a $\ttM $-algebra then $\bra
x,\alpha\c\lambda^\ttM _x,\alpha\c\rho^\ttM _x\ket$ is an $R$-$R$-bimodule.
\end{lem}

This provides the object map of a functor $\Psi\colon\M^\ttM \to \,_R\M_R$.
It is easy to show that if $\tau\colon\bra x,\alpha\ket\to\bra y,\beta\ket$ is
a $\ttM $-algebra morphism then it is also an $R$-$R$-bimodule morphism. This
defines a forgetful (monadic) functor $\Psi\colon\M^\ttM \to \,_R\M_R$. 
Moreover, the Eilenberg-Moore comparison functor $K\colon\C\to\M^\ttM $ 
satisfies $\Psi K=U$. Thus we can factorize the original functor $G$ in 3 steps
\beq
\begin{CD} \label{eq: fac}
\C@>K>>\M^\ttM@>\Psi>>_R\M_R@>\Gamma>>\M
\end{CD}
\eeq
Unfortunately, only the third functor is monoidal. In order to make it a
diagram in $\MonCat$ we need at least a monoidal structure on the
Eilenberg-Moore category $\M^\ttM$. 
This could be achieved under the stronger assumption that $U$ is strong
monoidal and has a left adjoint. Namely, \cite[Theorem 2.8]{Sz: JPAA} yields
the following result:
\begin{quote}
Let $U\colon\C\to\T$ be a strong monoidal right adjoint functor. Then 
\begin{itemize}
\item its monad $\ttT$ is opmonoidal, i.e., it is a monad on $\T$ in the
2-category $\OpmonCat$;
\item the category $\T^\ttT$ of $\ttT$-algebras has a unique monoidal structure
such that the Eilenberg-Moore forgetful functor $U^\ttT\colon\T^\ttT\to\T$ is
strict monoidal;
\item $U=U^\ttT K$ with a strong monoidal comparison functor
$K\colon\C\to\T^\ttT$.
\end{itemize}
\end{quote}
\begin{thm} \label{thm: Rmonfunc}
Assume that
\begin{enumerate}
\item $\bra\C,\b,e\ket$ is a monoidal category where $\C$ is complete, well
powered and has a small cogenerating set;
\item $\bra\M,\o,i\ket$ is a monoidal category where $\M$ has coequalizers and
$\o$ preserves them;
\item $\bra G,G_2,G_0\ket$ is an essentially strong monoidal
functor $\C\to\M$ with $G$ having a left adjoint.
\end{enumerate}
Then there is a monoid $R$ in $\M$ and an opmonoidal monad $\ttT$ on $_R\M_R$
such that
\beq \label{eq: mon fac}
\begin{CD}
\C@>K>>(_R\M_R)^\ttT\\
@VGVV @VV{U^\ttT}V\\
\M@<\Gamma<<_R\M_R
\end{CD}
\eeq
is commutative in $\MonCat$ where $K$ is the strong monoidal comparison
functor, $U^\ttT$ is the strict monoidal forgetful functor of $\ttT$-algebras
and $\Gamma$ is the monoidal forgetful functor of $R$-$R$-bimodules.
\end{thm}
\begin{proof}
By assumption (2) the canonical factorization $G=\Gamma U$ through $_R\M_R$
exists. Since $G$ is right adjoint, it preserves limits. But $\Gamma$, being
the forgetful functor of the monad $R\o\under\o R$ on $\M$, is monadic
therefore it creates limits. Therefore $U$ preserves limits, too. By the Special 
Adjoint Functor Theorem \cite[Corollary V. 8]{MacLane} assumption (1)
ensures that $U$ has a left adjoint. Now using assumption (3) the $U$ is
strong monoidal right adjoint, therefore, it factorizes monoidally through the
Eilenberg-Moore category of its bimonad $\ttT$ by the above quoted 
\cite[Theorem 2.8]{Sz: JPAA}. 
\end{proof}

The 3 step factorization found in the above Theorem describes what can be
expected in general for continuous monoidal functors. Of course, it would be
more interesting to replace the "abstract quantum groupoid" $\ttT$ with a
concrete bialgebroid, let us say. If $\M$ is the category of $k$-modules over a
commutative ring then a necessary and sufficient condition for $(_R\M_R)^\ttT$
to be the monoidal category $_A\M$ of modules over a bialgebroid $A$ was given
in \cite[Theorem 4.5]{Sz: JPAA}. The condition is very simple: the underlying
functor $T$ of the bimonad should have a right adjoint. Unfortunately, it is
difficult to find a condition on $G$ that guarantees a right adjoint for $U$.
In the next Section we will study a special case which allows to do so.

We note that even under the conditions of the Theorem the underlying monad of
$\ttT$ may not be $\ttM$, so the factorizations (\ref{eq: mon fac}) and
(\ref{eq: fac}) may be different. 


\section{Monoidal $\V$-functors to $\V$} \label{sec: bgd}

Now we turn to replace bimonads with bialgebroids. The key observation is
that every monoid $A$ in $\M$ determines a monad $A\o\under$ on $\M$ and -
under certain conditions on $\M$ - these monads are precisely the monads the
underlying endofunctor of which has a right adjoint. The Eilenberg-Moore
categories $\M^\ttT$ of such monads are precisely the categories $_A\M$ of
modules over $A$. This can be generalized to ($k$-linear) bimonads as follows.
\begin{thm}
\cite[Thm. 4.5]{Sz: JPAA} Let $\M=\M_k$ be the category of modules
over a commutative ring $k$,  $R$ be a $k$-algebra and $\ttT$ be a $k$-linear
bimonad on $_R\M_R$. Then there exists a bialgebroid $A$ in $\M$ over $R$ and
an isomorphism $\ttT\cong A\o\under$ if and only if $T\colon\M\to\M$ has a
right adjoint. Moreover, $T$ has a right adjoint if and only if
$U^\ttT\colon\M^\ttT\to\M$ has a right adjoint.  
\end{thm}
We want to give an analogue characterization of the long forgetful
functors $G\colon\,_A\M\to\M$ of bialgebroids. Clearly, these functors have
both left and right adjoints - the induction and coinduction functors - and the
factorization through $_R\M_R$ shows that they are essentially strong
monoidal. Thus we expect that these properties of a functor leads to the 
construction of a bialgebroid with Tannaka duality.

We let $\V$ denote either $_k\M$, $\Ab$ or $\Set$ and work with $\V$-monoidal
$\V$-categories and $\V$-functors between them that have $\V$-adjoints
\cite{Borceux II}. Due to the fact that $\V$ is not only symmetric monoidal
closed but its monoidal unit is a generator, we can work with the underlying
ordinary categories and consider $\V$-functors as a special class of ordinary
functors that preserve some extra structure that is encoded in the choice of
$\V$. The set of $\V$-natural transformations between $\V$-functors is the same
(under the 2-functor $\Hom_\V(i,\under)$) as the set of ordinary natural 
transformations.
So the complicated formalism of enriched categories can be avoided and
proofs of commutativity of diagrams in $\V$ can be done by \textit{elements}.
(An element of an object $v\in\V$ is an arrow $i\to v$ in $\V$.) In a
$\V$-category $\C$ we denote by $\C(a,b)\in\V$ its hom-objects and by
$\Hom_\C(a,b):=\Hom_\V(i,\C(a,b))$ its hom-sets. 
Otherwise it should be clear
from the context whether we speak about the $\V$-category $\C$ or its
underlying ordinary category. $\V$ itself is a $\V$-category with $\V(x,y)$
being the internal hom object $\hom(x,y)$ which is defined by
\[
\Hom_\V(v\o x,y)\cong\Hom_\V(v,\hom(x,y))\,.
\]
Our base category $\V$ is also complete which allows to define the Takeuchi
$\x_R$-product as a pullback of equalizers. Existence of coequalizers is
also needed in $\V$ if we want to apply Theorem \ref{thm: monfunc}.
Readers interested in more general enrichments than the three cases mentioned
above can take for $\V$ any complete category with coequalizers which is
endowed with a symmetric closed monoidal structure $\bra\V,\o,i\ket$ such that
$i$ a projective generator. 

We restrict ourselves to study functors with target category being the base
category, i.e., $G$ is a $\V$-functor to $\V$. In this
situation the existence of a left $\V$-adjoint $F\dashv G$ implies that
$G\colon\C\to\V$ is representable: 
\beq
Ga\iso\hom(i,Ga)\iso\C(Fi,a)
\eeq
for $a\in\C$. Thus without loss of generality we may assume that $G$ is a
hom-functor $G=\C(g,\under)$. If such a $G$ has a monoidal structure then - by
the Yoneda Lemma - a comonoid structure $\bra g,\gamma,\pi\ket$ on $g$
arises via
\begin{alignat}{2}
G\colon&\C\to\V&\qquad G&=\C(g,\under)\\
G_{a,b}\colon& Ga\o Gb\to G(a\b b)&\quad
G_{a,b}(\alpha\o\beta)&=(\alpha\b\beta)\c\gamma\\
G_0\colon& i\to Ge&\qquad G_0&=\pi
\end{alignat}
where in the last equation one can recognize the usual identification of
elements of a hom-object of a $\V$-category with arrows in $\V$.

The monoid $R$ in the canonical factorization $G=\Gamma U$ through
$_R\V_R$ is nothing but the convolution monoid $\C(g,e)$ with multiplication 
\beq
\rho\ast \rho'=\l_e\c(\rho\o \rho')\c\gamma\,,\qquad\rho,\rho'\in R
\eeq
and unit $\pi$.

We note also that every object $Ga=\C(g,a)\in\V$ has a natural right $A$-module
structure where $A:=\C(g,g)$ is the endomorphism monoid in $\V$. So $G$
factorizes through the forgetful functor $\V_A\to\V$. This is compatible with
the canonical factorization $G=\Gamma U$ because we have monoid morphisms
\begin{alignat}{2}
s\colon& R\to A&\qquad
s(\rho)&:=\r_g\c (g\b \rho)\c \gamma\\
t\colon& R\op\to A&\qquad
t(\rho)&:=\l_g\c(\rho\b g)\c\gamma
\end{alignat}
Corresponding to the diagram $i\to R\op\o R\rarr{t\o s}A$ in $\Mon\V$ there is
the diagram $\V_A\to\,_R\V_R\to\V$ in $\V$-$\Cat$.

\begin{thm} \label{thm: RVmonfunc}
Denoting by $\V$ either $\M_k$, or $\Ab$ or $\Set$
let $\C$ be a $\V$-monoidal $\V$-category and $G\colon\C\to\V$ be a
representable essentially strong monoidal $\V$-functor. 
Then there exists a
monoid $R$ in $\V$, a right bialgebroid $A$ in $\V$ over $R$, a strong monoidal
$K\colon\C\to\V_A$ and a monoidal natural isomorphism $G\iso \Gamma^R\,\Phi^A\,
K$, 
\beq 
\begin{CD}
\C@>K>>\V_A\\
@V{G}V{\hskip 0.7truecm\cong}V @VV{\Phi^A}V\\
\V@<\Gamma^R<<_R\V_R
\end{CD}
\eeq
where $\Phi^A$ denotes the strict monoidal forgetful functor of the
bialgebroid $A$.
\end{thm}
\begin{proof}
As we have explained above, there exists a representing comonoid $\bra
g,\gamma,\pi\ket$ in $\C$ and a monoidal natural isomorphism
$G\iso\C(g,\under)$. 
As a $\V$-functor, $\C(g,\under)$ factorizes as
\beq
\begin{CD}
\C@>K>>\V_A@>\Phi>>_R\V_R@>\Gamma>>\V
\end{CD}
\eeq
where $A=\C(g,g)$ is the endomorphism monoid, $R=\C(g,e)$ is the
convolution monoid and - identifying $_R\V_R$ with $\V_{R^e}$ where $R^e=R\op\o
R$ - $\Phi$ is the $\V$-functor that forgets along $t\o s\colon R^e\to A$. By
the very definition of $R$ we see that the canonical factorization $\Gamma U$
of $\C(g,\under)$ leads to $U=\Phi K$ which is strong monoidal, hence
opmonoidal. Since opmonoidal functors map comonoids to comonoids, the triple 
\begin{align} 
C&=Ug\equiv\bra A,R\o A\iso A\o R\rarr{A\o t}A\o
A\rarr{\mu}A,A\o R\rarr{A\o s}A\o A\rarr{\mu}A\ket\\
\delta&=\left(C\rarr{U\gamma}U(g\b g)\rarr{U_{g,g}^{-1}}C\o_R C\right)\\
\eps&=\left(C\rarr{U\pi}Ue\rarr{U_0^{-1}}\,_RR_R\right)
\end{align}
is a comonoid in $_R\V_R$.
This comonoid allows to define a monoidal structure on $\V_A$ such that
$\Phi$ becomes strict monoidal. As a matter of fact, let $X$ and
$Y$ be right $A$-modules in $\V$ and define the $A$-module
$X\b_A Y$ as the object $X\o_R Y$ with $A$-action
\beq
(x\o_R y)\ract\alpha:=(x\ract\alpha\1)\o_R (y\ract\alpha\2)\,,\qquad 
x\in X,y\in Y,\alpha\in A\,. 
\eeq
But in order for this action to be well-defined we have to show that
$\delta(\alpha)$, $\alpha\in A$, belong to the subbimodule $C\x_R C\subset C\o_R
C$ defined as the intersection of the small set of equalizers
$\{e_\rho\,|\,\rho\in R\}$ where 
\beq
C\x_\rho C\ \ \stackrel{e_\rho}{\rightarrowtail}\ \ C\o_R C
\begin{gathered}
Us(\rho)\o_R C\\
\begin{picture}(70,8)
\put(10,2){\vector(1,0){50}}
\put(10,5){\vector(1,0){50}}
\end{picture}\\
C\o_R Ut(\rho)
\end{gathered}
C\o_R C
\eeq
Here $Us(\rho)=s(\rho)\c\under$ and $Ut(\rho)=t(\rho)\c\under$ 
so they are $R$-$R$-endomorphisms of $C$. 
Naturality of $U_2$, the definition of $\delta$ and some elementary monoidal
calculus yields the following identity in $_R\V_R$, 
\begin{align*}
U_{g,g}\c (Us(\rho)\o_R Ug)\c\delta&=U(s(\rho)\b g)\c U_{g,g}\c\delta\\
&=U(s(\rho)\b g)\c U\gamma\\
&=U((\r_g\b g)\c((g\b \rho)\b g)\c (\gamma\b g)\c\gamma)\\
&=U((g\b\l_g)\c(g\b(\rho\b g))\c(g\b\gamma)\c\gamma)\\
&=U(g\b t(\rho))\c U\gamma\\
&=U_{g,g}\c(Ug\o_R Ut(\rho))\c\delta\ .
\end{align*}
Since $U_{g,g}$ is invertible, $\delta$ restricts to a map $A\to
A\x_R A$. 
Notice that $A\x_R A$ - as an object in $\V$ - inherits a monoid structure
from that of $A\o A$ and then $\delta$ becomes a morphism of monoids. As a
matter of fact,
\begin{align*}
U_{g,g}(\delta(\alpha)\delta(\beta))&=(\alpha\1\b\alpha\2)\c\gamma\c\beta\\
&=\gamma\c\alpha\c\beta\\
&=U_{g,g}(\delta(\alpha\c\beta))
\end{align*}
for $\alpha,\beta\in A$. Unitality of $\delta$ is obvious.
This finishes the definition of $\b_A$ and the associativity coherence
isomorphism $\a$ of $\V_A$ as a lift of $\a$ of $_R\V_R$. Whether it has a
unit object depends on the counit properties. 
\begin{align*}
\eps(s(\eps(\alpha))\c\beta)&=\pi\c s(\eps(\alpha))\c\beta\\
&=\pi\c\r_g\c(g\b\eps(\alpha))\c\gamma\c\beta\\
&=\r_e\c(\pi\b\eps(\alpha))\c\gamma\c\beta=\eps(\alpha\c\beta)=\pi\c\alpha\c\beta\\
&=\eps(\alpha\c\beta)
\end{align*}
and similarly, $\eps(t(\eps(\alpha))\c\beta)=\eps(\alpha\c\beta)$ holds for
all $\alpha,\beta\in A$. Unitality of $\eps$ is obvious since $\eps(g)=\pi$ is
the unit element of $R$. Then the monoidal unit $E$ of $\V_A$ becomes 
$\bra R,R\o A\rarr{s\o A}A\o A\rarr{\mu}A\rarr{\eps}R\ket$ and the unit
coherence isomorphisms $\l$ and $\r$ of $_R\V_R$ lift to become
the $\l$ and $\r$ of $\V_A$, respectively.

This finishes the proof of that $A$ is a bialgebroid together with
the construction of a monoidal structure on $\V_A$ such that $\Phi$ is strict
monoidal. Since $U$ is strong and $\Phi$ reflects isomorphisms, it follows
that $K$ is strong monoidal.

\end{proof}


\section{Characterizing long forgetful functors of bialgebroids} 
\label{sec: G^A}

Recall that for every closed monoidal category $\V$ there is a monoidal
equivalence (of ordinary categories)
\beq
\V\simeq \mathsf{L}\text{-}\V\text{-}\Fun(\V,\V)\,,\qquad v\mapsto \,\under\o v
\eeq
of $\V$ with the strict monoidal category of left adjoint $\V$-functors
$\V\to\V$. 

Since the $\V$-monads on $\V$ the underlying $\V$-functors of which
have right $\V$-adjoints are precisely the monoids in 
$\mathsf{L}\text{-}\V\text{-}\Fun(\V,\V)$, they are mapped by the above
equivalence into the monoids in $\V$. 

Let $\V\text{-}\Cat/\V$ denote the 2-category with objects the $\V$-functors
to $\V$ and with arrows $F\to G$ the $\V$-functors $K$ with an isomorphism
$F\cong GK$. One defines $\V\text{-}\MonCat/\V$ similarly. 

For a right bialgebroid $A$ in $\V$ we denote by $G^A$ the 
essentially strong monoidal forgetting $\V$-functor $\V_A\to\V$.
$G^A$ is called the long forgetful functor of $A$.
\begin{thm} \label{thm: G^A}
A $\V$-functor $G\colon\C\to\V$ is equivalent in $\V\text{-}\Cat/\V$ 
to the long forgetful functor $G^A$ of a bialgebroid in $\V$ iff 
\begin{itemize} 
\item $G$ is monadic,
\item $G$ has a right adjoint and
\item there is an essentially strong monoidal structure $\bra G,G_2,G_0\ket$
on $G$.
\end{itemize}
In this case $G\simeq G^A$ in $\V\text{-}\MonCat/\V$.
\end{thm}
\begin{proof}
Necessity: Since the object in $\V$ underlying the monoidal unit of $\V_A$ is
just the object underlying the base monoid $R$, the canonical factorization of
$G^A$ is $G^A=\Gamma^R\Phi^A$, in the notation of Theorem \ref{thm: RVmonfunc}.
Thus $G^A$ is essentially strong.
As every forgetful functor of a monoid, $G^A$ is monadic.
It has a right adjoint, namely, the coinduction functor sending the object $x$
of $\V$ to the object $\hom(A,x)$ endowed with right $A$-action
\beq
\eps_x\colon\hom(A,x)\o A\to x
\eeq
provided by the $x$-component of the counit of the adjunction $\under\o
A\dashv\hom(A,\under)$.

Sufficiency: Since $G$ is monadic, it has a left adjoint $F$ and $G=U^\ttM J$
where $\ttM$ is the monad with underlying functor $GF$, $U^\ttM$ is its
forgetful functor and $J\colon\C\to \V^\ttM$ is an equivalence. Let $H$ be the
right adjoint of $G$. Then $GH$ is a right adjoint of $GF$, therefore, by the
above remark, it is isomorphic to the monad $\under\o B$ associated to a monoid
$B$ in $\V$. Therefore $\V^\ttM$ is isomorphic to the category $\V_B$ of right
$B$-modules and the decomposition $G=U^\ttM J$ can be replaced with $G=G^B L$
where $L\colon\C\to\V_B$ is an equivalence. Using the latter equivalence
$G^B$ is given an essentially strong monoidal structure and it has a left
adjoint. Therefore Theorem \ref{thm: RVmonfunc} provides a monoidal isomorphism
$G^B\cong\Gamma\Phi K=G^A K$ where the right bialgebroid $A$ has underlying
monoid the endomorphism monoid of the representing object $g=B_B$ of $G^B$.
Clearly, $A=\End(B_B)\cong B$. This proves that $K\colon \V_B\to\V_A$,
$X_B\mapsto \Hom_B(\,_AB_B,X_B)$, is an isomorphism of (monoidal) categories.
Since in the monoidal isomorphism $G\cong G^A KL$ the functor $KL\colon
\C\to\V_A$ is a monoidal equivalence, it defines an equivalence $G\iso G^A$ in
$\V\text{-}\MonCat/\V$.  
\end{proof}

\section{The forgetful functors of weak bialgebras} \label{sec: WBA}

In this section the base category $\V$ is the category $\M_k$ of modules
over a commutative ring. So we switch to the convention of
writing capital Roman letters for objects and corresponding small case
letters for their elements.
Weak bialgebras over $k$ are not only special bialgebroids in $\M_k$ but are
equipped with some more structure, as well. The extra structure can be
recognized in two places: (1) in the difference between a weak bialgebra counit
$\epsilon\colon B\to k$ and a bialgebroid counit $\eps\colon A\to R$ and (2)
in the nontrivial element $\cop(1)\in A\o A$ which is closely related to a
separability idempotent of the separable algebra $R$. In other words, a weak
bialgebra $A$ is a bialgebroid over a separable Frobenius algebra $R$ together
with a Frobenius structure $\bra\epsilon,\sum_i e_i\o f_i\ket$ in which $\sum_i
e_if_i=1$ (see \cite{Sz: Strbg} for more details). We call $\bra
R,\epsilon,\sum_i e_i\o f_i\ket$ a separable Frobenius structure. 

Accordingly, the forgetful functor $G^A\colon\M_A\to\M_k$ of a bialgebroid has
more structure than just a monadic essentially strong monoidal functor with
right adjoint. It is equipped also with an opmonoidal structure. For
concreteness let $\bra A,\cop,\epsilon\ket$ be a weak bialgebra \cite{BNSz}
and let $R$ be identified with the canonical right subalgebra
$A^R=\{1\1\epsilon(a1\2)\,|\,a\in A\}$. Then $A$ becomes a right bialgebroid
over $R$ with
\begin{alignat}{2}
s\colon& R\to A\,,&\qquad r&\mapsto r\\
t\colon& R\op\to A\,,&\qquad r&\mapsto \epsilon(r1\1)1\2\\
\delta\colon&A\to A\o_R A\,,&\qquad \delta&:=\tau\c\cop\\
\eps\colon&A\to R\,,&\qquad a&\mapsto 1\1\epsilon(a1\2)
\end{alignat}
where $\tau\colon A\o A\to A\o_R A$ is the canonical epimorphism (cf.
\cite[Lemma 1.1]{Sz: Strbg}). 

Let $G$ denote the long forgetful functor $G^A$ of the bialgebroid $A$ (see
Section \ref{sec: bgd}). Then the monoidal structure 
\begin{alignat}{2}
G_{X,Y}\colon& GX\o GY\to G(X\b_A Y)\,,&\qquad x\o y&\mapsto x\o_R y\\
G_0\colon& k\to GE\,,&\qquad \kappa&\mapsto\kappa\cdot 1_R
\end{alignat}
is built of the canonical epimorphisms $X\o Y\to X\o_R Y$ for
$R$-$R$-bimodules and of the unit $k\to R$ of the algebra $R$. Due to the
presence of the separable Frobenius structure it has an opmonoidal counterpart 
\begin{alignat}{2}
G^{X,Y}\colon& G(X\b_A Y)\to GX\o GY\,,&\qquad x\o_R y&\mapsto \sum_i x\cdot
e_i\o f_i\cdot y\\
G^0\colon& GE\to k\,,&\qquad r&\mapsto \epsilon(r)
\end{alignat}
So we have a monoidal $\bra G,G_2,G_0\ket$ and an opmonoidal functor $\bra G,
G^2,G^0\ket$ with the same underlying functor $G$. This involves that $G_2$
with $G_0$ obey one hexagonal and two square diagrams and $G^2$ together with
$G^0$ obey three analogous, but oppositely oriented, diagrams. In addition,
there are compatibility conditions between the monoidal and opmonoidal
structures that are not of the bialgebra type but rather of the Frobenius
algebra type \cite{Abrams}. Namely, 
\begin{align} 
(G_{X,Y}\o GZ)\c\a_{GX,GY,GZ}\c(GX\o G^{Y,Z})
&=G^{X\b Y,Z}\c G\a_{X,Y,Z}\c G_{X,Y\b Z} \label{eq: frob1}\\
(GX\o G_{Y,Z})\c\a^{-1}_{GX,GY,GZ}\c(G^{X,Y}\o GZ)
&=G^{X,Y\b Z}\c G\a^{-1}_{X,Y,Z}\c G_{X\b Y,Z} \label{eq: frob2}
\end{align}
for all $A$-modules $X$,$Y$, $Z$ and where $\b$ stands for the monoidal
product $\b_A$ of $A$-modules. At last but not least, $G_2$ is split epi and
the splitting map is just $G^2$, i.e.,
\beq \label{eq: sep}
G_{X,Y}\c G^{X,Y}\ =\ G(X\b Y)\,,\qquad X,Y\in\M_A\,.
\eeq

Now we summarize these experiences in a
\begin{defi}
Let $\bra\C,\b,E\ket$ and $\bra\M,\o,I\ket$ be monoidal categories.
A separable Frobenius structure on a functor $G\colon\C\to\M$ consists of
natural transformations $G_2$, $G^2$ and arrows $G_0$, $G^0$ such that
\begin{enumerate}
\item $\bra G,G_2,G_0\ket$ is a monoidal functor,
\item $\bra G,G^2,G^0\ket$ is an opmonoidal functor,
\item the Frobenius conditions (\ref{eq: frob1}), (\ref{eq: frob2}) and the
separability condition (\ref{eq: sep}) hold.
\end{enumerate}
\end{defi}
Of course, the functor itself may be neither separable nor Frobenius. Only 
its monoidal-opmonoidal structure is restricted in the above Definition. 
\begin{lem} \label{lem: split}
If $\bra G,G_2,G_0,G^2,G^0\ket$ is a separable Frobenius structure on the
functor $G\colon\C\to\M$ then (\ref{exact}) is a split coequalizer in $\M$ for
all pairs of objects in $\C$. In particular, $\bra G,G_2, G_0\ket$ is
essentially strong monoidal.  
\end{lem}
\begin{proof}
For each pair $X,Y$ of objects in $\C$ we can define arrows in $\M$ by
\begin{align*}
\partial_0&:=(GX\o G\l_Y)\c(GX\o G_{E,Y})\\
\partial_1&:=(G\r_X\o GY)\c(G_{X,E}\o GY)\c\a_{GX,GE,GY}\\
\gamma&:=G_{X,Y}\\
\sigma&:=G^{X,Y}\\
\tau&:=(GX\o G^{E,Y})\c(GX\o G\l^{-1}_Y)
\end{align*}
Then an elementary calculation gives
\begin{alignat*}{2}
&\text{hexagon for $G_2$}\qquad
&\Rightarrow\qquad\gamma\c\partial_0&=\gamma\c\partial_1\\
&\text{separability (\ref{eq: sep})}\qquad
&\Rightarrow\qquad \gamma\c\sigma&=1\\
&\text{separability (\ref{eq: sep})}\qquad
&\Rightarrow\qquad \partial_0\c\tau&=1\\
&\text{Frobenius (\ref{eq: frob1})}\qquad
&\Rightarrow\qquad\partial_1\c\tau&=\sigma\c\gamma  
\end{alignat*}
which means precisely that $\gamma$ is a split coequalizer, hence a
coequalizer, of $\partial_0$ and $\partial_1$ \cite{MacLane}. Thus $G$ is
essentially strong monoidal.
\end{proof}

\begin{lem} \label{lem: sep Frob}
Let $\bra G,G_2,G_0,G^2,G^0\ket$ be a separable Frobenius structure on the
functor $G\colon\C\to\M$ between monoidal categories. Then the image $GE$ of
the unit object of $\C$ is equipped with a separable Frobenius structure in
$\M$.  
\end{lem}
\begin{proof}
$GE$ gets a monoid structure as the image of $E$ by $\bra G,G_2,G_0\ket$ as we
explained in Lemma \ref{lem: mon to mon}. Dually, the $\bra G,G^2,G^0\ket$ maps
the comonoid $E$ into a comonoid with underlying object $GE$. So we obtain
\begin{align}
R&=GE\\
\mu&=G\l_E\c G_{E,E}\ :\ R\o R\to R\\
\eta&=G_0\ :\ I\to R\\
\sigma&=G^{E,E}\c G\l_E^{-1}\ :\ R\to R\o R\\
\psi&=G^0\ :\ R\to I
\end{align}
such that $\bra R,\mu,\eta\ket$ is a monoid and $\bra R,\sigma,\psi\ket$ is a
comonoid in $\M$. Now (\ref{eq: frob1}), (\ref{eq: frob2}) and 
(\ref{eq: sep}) imply that these two structures on $GE$ are compatible 
in the sense of satisfying
\begin{align}
(\mu\o R)\c \a_{R,R,R}\c(R\o\sigma)&=\sigma\c\mu\\
(R\o\mu)\c\a_{R,R,R}^{-1}\c(\sigma\o R)&=\sigma\c\mu\\
\mu\c\sigma&=R
\end{align}
which are the defining relations of a separable Frobenius structure on $R$. 
\end{proof}
We note that for the familiar categories, $\M=\M_k$, the $\sigma$ is uniquely
determined by $\sigma(1_R)$ which in turn is uniquely determined by $\psi$. 
If, moreover, $k$ is a field then a separable Frobenius
algebra, i.e., an algebra having a separable Frobenius structure (also called
an index one Frobenius algebra), is nothing but a separable $k$-algebra
\cite{KSz}. 

\begin{lem} \label{lem: sep bim}
Let $\bra R,\mu,\eta,\sigma,\psi\ket$ be a separable Fobenius structure in 
$\M_k$. Then the monoidal forgetful functor of bimodules
$\Gamma=\Gamma^R\colon\,_R\M_R\to\M$ has the following extension to a separable
Frobenius structure:
\begin{alignat}{2}
\Gamma^{X,Y}\colon& X\o_R Y\to X \o Y &\qquad
x\o_R y&\mapsto \sum_i x\cdot e_i\o f_i\cdot y_i\\
\Gamma^0\colon& R\to k&\qquad r&\mapsto \psi(r)
\end{alignat}
where $\sum_i e_i\o f_i=\sigma(1_R)$. 
\end{lem}
\begin{proof}
This is left for an exercise.
\end{proof}

After the above preparations Theorem \ref{thm: G^A} has the following
\begin{cor}
A $k$-linear functor $G\colon \C\to \M_k$ - as an object in 
$k\text{-}\Cat/\M_k$ - is equivalent to the long forgetful functor of a weak
bialgebra iff 
\begin{itemize}
\item $G$ is monadic,
\item $G$ has a right adjoint and
\item there is a separable Frobenius structure $\bra G,G_2,G_0,G^2,G^0\ket$ on
$G$. 
\end{itemize}

\end{cor}
\begin{proof}
$G$ is essentially strong monoidal by Lemma \ref{lem: split} therefore
Theorem \ref{thm: G^A} provides a right bialgebroid 
$\bra A,R,s,t,\delta,\eps\ket$ and a monoidal equivalence $G\simeq G^A$.
It remains to show that the data on $A$ can be extended to the data of a
weak bialgebra. (This extra structure is encoded neither in the monoidal
category $\M_A$ nor in the monoidal functor $G^A$.) A monoidal equivalence is
the same thing as an opmonoidal equivalence, so we can use $G\simeq G^A$ to
pass the whole separable Frobenius structure of $G$ to $G^A$. Now Lemma
\ref{lem: sep Frob} implies that $R=G^AE$, the base of $A$, is given a
separable Frobenius algebra structure $\bra R,\mu,\eta,\sigma,\psi\ket$. 
Then weak bialgebra comultiplication and counit can be introduced by
\begin{alignat}{2}
\cop&:=\Gamma^{A,A}\c\delta&\qquad\cop(a)&=\sum_i a\1 s(e_i)\o a\2 t(f_i)\,,\\
\epsilon&:=\Gamma^0\c\eps&\qquad\epsilon(a)&=\psi(\eps(a))\,.
\end{alignat}
For the proof of that $\bra A,\cop,\epsilon\ket$ is a weak bialgebra 
we refer to the proof of \cite[Proposition 7.4]{KSz} where this has been
done for $R$ a separable $k$-algebra over a field. However, after providing a
separable Frobenius structure on $R$ that proof applies here.  
\end{proof}

\end{document}